\documentclass[letterpaper, 10 pt, conference]{ieeeconf}  

\IEEEoverridecommandlockouts                              

\usepackage{graphics} 
\usepackage{epsfig} 
\usepackage{times} 
\usepackage{mathtools}
\usepackage{amsmath} 
\usepackage{amssymb}  
\usepackage{amsthm}
\usepackage{color}

\usepackage{tikz}
\usetikzlibrary{shapes,arrows}
\usepackage{pgfplots} 
\usepackage{pgfgantt}
\usepackage{pdflscape}
\usepackage{nicefrac}
\usepackage{prettyref}
\pgfplotsset{compat=1.5}
\pgfplotsset{plot coordinates/math parser=false}
\newlength\fwidth
\usepackage{grffile}
\usetikzlibrary{plotmarks}
\usetikzlibrary{arrows.meta}
\usepgfplotslibrary{patchplots}
\usepackage{units}

\newtheorem{definition}{Definition}
\newtheorem{remark}{Remark}
\newtheorem{theorem}{Theorem}
\newtheorem{lemma}{Lemma}
\newtheorem{problem}{Problem}
\newtheorem{proposition}{Proposition}
\newtheorem{assumption}{Assumption}

\newcommand{\defeq}{\vcentcolon=}
\newcommand{\T}{\scriptscriptstyle\top}       
\newcommand{\mathmin}{\operatorname*{min}}

\newcommand{\mathst}{\text{s.t.}}

\definecolor{myred}{RGB}{233,72,73}%
\definecolor{mygreen}{RGB}{113,191,110}%
\definecolor{myblue}{RGB}{93,147,191}%
\definecolor{mydarkblue}{RGB}{57,101,181}%
\definecolor{mycyan}{rgb}{ 0.05, 0.80, 0.75}%
\definecolor{mygray}{RGB}{163,163,163}%
\definecolor{myorange}{RGB}{233,180,73}%

\long\def\comment#1{}
\renewenvironment{proof}{{\bfseries Proof:~}}{}

\usepackage{eso-pic}
\AddToShipoutPictureBG*{%
\AtPageUpperLeft{%
\setlength\unitlength{1in}%
\hspace*{\dimexpr0.5\paperwidth\relax}
\makebox(0,-21)[c]{
\footnotesize
\begin{tabular}{c c}
© 2023 IEEE. Personal use of this material is permitted. Permission from
IEEE must be obtained for all other uses, in any \\
current or future media, including reprinting/republishing
this material for advertising or promotional purposes, creating new \\
collective works, for resale or redistribution to servers or lists,
or reuse of any copyrighted component of this work in other works.
\end{tabular}}}}

\usepackage[super]{nth}
\AddToShipoutPictureBG*{%
\AtPageUpperLeft{%
\setlength\unitlength{1in}%
\hspace*{\dimexpr0.5\paperwidth\relax}
\makebox(0,-1.4)[c]{
\begin{tabular}{c c}
To appear in the Proceedings of the 6\nth{2} IEEE Conference on Decision and Control, Singapore, 2023. \\
Uploaded to arXiv on \today. 
\end{tabular}}}}

\title{\LARGE \bf
 Discrete-time Control Barrier Functions for Guaranteed Recursive Feasibility in Nonlinear MPC: An Application to Lane Merging
}

\author{Alexander Katriniok$^{1}$, Erfan Shakhesi$^{1}$ and W.P.M.H. (Maurice) Heemels$^{1}$
\thanks{$^{1}$Alexander Katriniok, Erfan Shakhesi and Maurice Heemels are with the Control Systems Technology section, Mechanical Engineering, Eindhoven University of Technology, 
The Netherlands,
        {\tt\small a.katriniok@tue.nl, e.shakhesi@student.tue.nl, m.heemels@tue.nl}}%
}

\newrefformat{sec}{\mbox{Section \ref{#1}}}
\newrefformat{tab}{\mbox{Table \ref{#1}}}
\newrefformat{fig}{\mbox{Fig. \ref{#1}}}
\newrefformat{rem}{\mbox{Remark \ref{#1}}}
\newrefformat{alg}{\mbox{Algorithm \ref{#1}}}
\newrefformat{ass}{\mbox{Assumption \ref{#1}}}
\newrefformat{def}{\mbox{Definition \ref{#1}}}
\newrefformat{prop}{\mbox{Proposition \ref{#1}}}
\newrefformat{prob}{\mbox{Problem \ref{#1}}}
\newrefformat{thm}{\mbox{Theorem \ref{#1}}}

\begin{document}

\maketitle
\thispagestyle{empty}
\pagestyle{empty}

\begin{abstract}
In this paper, we present conditions under which the terminal ingredients, defined by discrete-time control barrier function (DTCBF) certificates, guarantee recursive feasibility in nonlinear MPC. Further, we introduce the notion of \mbox{quasi-DTCBF} (qDTCBF) certificates. Compared to DTCBFs, qDTCBF conditions can be satisfied with tighter control input bounds, which is highly advantageous if only limited actuation is possible. Both certificates encourage an earlier reaction of the control system and result in a lower cumulative MPC cost. The methodology is applied to a lane merging problem in automated driving, in which DTCBF and qDTCBF certificates subject to input constraints form the terminal ingredients to guarantee recursive feasibility of the nonlinear MPC scheme. A simulation study demonstrates the efficacy of the concept. 
\end{abstract}

\section{INTRODUCTION}
\label{sec:introduction}


Rigorous safety guarantees are fundamentally important in the control of safety-critical systems, which can be found in many domains, especially those where humans are involved, including medical systems, critical infrastructures and autonomous vehicles (AVs). In this paper a special focus will be on AVs, which have to avoid collisions with other road users (agents) while accomplishing their driving task \cite{Yu2021a}. 

In AV motion planning, optimization-based methods such as model predictive control (MPC) are an appealing choice to explicitly accommodate constraints or exploit anticipated trajectories of other agents. The design of 
suitable terminal ingredients 
allows to guarantee recursive feasibility \cite{Rawlings2022a} (and thus safety). However, it is oftentimes challenging to solve the underlying, mostly nonconvex optimization problems in real-time. In recent years, control barrier functions (CBFs) have gained increasing attention as an alternative method to rigorously guarantee constraint satisfaction in an optimal control framework \cite{Ames2019a,Katriniok2022a}. Although computationally very efficient, disadvantages include rather instantaneous (harsh) control actions compared to MPC, the difficulty to utilize predictions and their proneness to gridlocks \cite{Ames2019a,Katriniok2022a,Wang2017a}. That said, the combination of the two methods appears to be a promising research direction. 

CBFs have originally been introduced for continuous-time systems \cite{Ames2014a}. In controller synthesis, linear CBF constraints (or certificates) are imposed in a quadratic programming (QP) problem to determine the least deviation from a nominal control input to render a safe set controlled invariant. An extension towards discrete-time systems is discussed in \cite{Agrawal2017a}. The authors in \cite{Xiong2022a} outline further extensions to accommodate higher relative degree and to relax input constraints through penalty functions for enhanced feasibility. Initial work to integrate MPC with discrete-time CBFs (DTCBFs) is proposed in \cite{Zeng2021a}, where DTCBF-inspired certificates (\textit{inspired} means that contrary to the definition in \cite{Agrawal2017a}, the safe set is not controlled invariant) are imposed over the entire prediction horizon to encourage an earlier controller reaction to obstacles and reduce the cumulative actuation cost. 
To enhance feasibility, 
\cite{Zeng2021b} introduces the decay rate of these certificates as additional decision variable. 
Recursive feasibility, though, is not addressed in \cite{Zeng2021a,Zeng2021b}.
A recursively feasible robust tube-based MPC scheme for nonlinear sampled-data systems, applying a DTCBF certificate at the terminal stage, is presented in \cite{Schillinger2021a}. In this setting, robust tubes accommodate the error of numerical integration. The control scheme, though, requires the initial state to belong to a (tightened) safe set and joint feasibility of the DTCBF certificate with input constraints is not guaranteed by construction, but is part of the assumptions. An MPC-based safety filter, recovering the control system from infeasible configurations through a soft-constrained terminal DTCBF \mbox{certificate, is outlined in \cite{Wabersich2022a}.}


\vspace*{1mm}\noindent
\textit{Main Contribution:} 
We propose to utilize DTCBF and novel quasi-DTCBF certificates to guarantee recursive feasibility in nonlinear MPC (NMPC), and apply this methodology to a lane merging application. This leads to the following contributions compared to existing works in literature: 
\begin{itemize} 
\item We convey conditions under which the terminal ingredients, defined either by \mbox{\textsc{i}) DTCBF} or \mbox{\textsc{ii}) quasi-DTCBF} certificates (weaker condition), guarantee recursive feasibility of the NMPC 
scheme. This aspect is formally not covered by \cite{Zeng2021a,Zeng2021b}.
\item Terminal certificates \textsc{i}) and \textsc{ii}) introduce an additional degree of freedom to reduce the cumulative tracking, actuation and stage cost through an earlier controller reaction. Compared to \textsc{i}), \textsc{ii}) is feasible for tighter input bounds, which is beneficial in case of limited actuation.
\item Unlike \cite{Zeng2021a,Zeng2021b,Schillinger2021a}, the initial state does not have to be within the safe set associated with these certificates. In fact, we allow the state to safely transition into the safe set over the horizon to reduce conservatism. 
\item For lane merging, recursive feasibility has been addressed in \cite{Bali2018,Geurts2023a} following a \textit{classical} controlled invariant set approach in a mixed-integer MPC setting. We formulate the problem without integer variables and apply (quasi-)DTCBF certificates, while guaranteeing joint feasibility of multiple certificates subject to input constraints, which is not addressed in \cite{Schillinger2021a,Wabersich2022a}.
\end{itemize} 
\noindent
\textit{Outline:} 
\prettyref{sec:problem} defines the control problem, which is solved in an NMPC setting as described in \prettyref{sec:NMPC}. Then, \prettyref{sec:DTCBF} introduces the theory on (quasi-)DTCBFs to guarantee recursive feasibility in NMPC,   
while \prettyref{sec:NMPCwCBF} presents its application to lane merging. Simulation results are then discussed in \prettyref{sec:results}. Finally, \prettyref{sec:conclusion} provides a conclusion and an outlook to future work.

\section{LANE MERGING PROBLEM}
\label{sec:problem}

\subsection{Problem Description}
\label{sec:problem_description}

The lane merging problem to solve (see \prettyref{fig:problem_description_sketch}), involves two agents $i \in \mathcal{A} \defeq \{1,2\}$ and can be stated as follows.

\begin{problem}[Lane Merging] \normalfont
\label{prob:problem_description_problemDef}
Agent 1 and Agent 2 drive in two separate lanes, each agent with its own reference velocity to track. Given \textit{a priori} known paths $\mathcal{P}_1$ and $\mathcal{P}_2$, \mbox{Agent 1} enters the lane of Agent 2 at the lane change point while the merging point constitutes the end of the merging maneuver. A recursively feasible centralized control scheme shall guarantee a safe distance between agents at all times, and respect the bounds on the agents' velocity and their control inputs.
\end{problem}

\begin{assumption} \normalfont 
\label{ass:problem_description_assumption1} 
To limit 
complexity of \prettyref{prob:problem_description_problemDef}, we rely on the following assumptions:
\begin{enumerate}
\item The agents' path $\mathcal{P}_i$ is \textit{a priori} known and fixed.
\item Agents follow their path without lateral displacement.
\item Only longitudinal dynamics are controlled.
\item Agents only drive in the forward direction.
\end{enumerate}
\end{assumption}

\subsection{Modeling}
\label{sec:problem_modeling}

\begin{figure}[t!]
	\centering
	\def\svgwidth{80mm}
	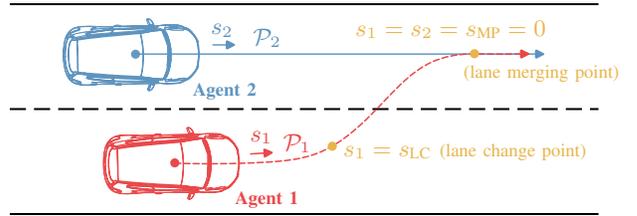	
	\vspace*{-2mm}
	\caption{Lane merging problem with 2 agents. Agent 1 enters the lane of Agent 2 at the lane change (LC) point and finishes the maneuver at the merging point (MP). At MP, 
    it holds $s_1 = s_2 = s_{\text{MP}}\defeq 0$.} 
	\vspace*{-5mm}
\label{fig:problem_description_sketch}
\end{figure}

With \prettyref{ass:problem_description_assumption1}, the dynamics of Agent $i$ along the \textit{a priori} known path $\mathcal{P}_i$ is stated as a double integrator model  
\begin{align} \label{eq:problem_modeling_sysContTime}
	\dot{s}_i = v_i, ~~ \dot{v}_i = a_{i},~~~ i \in \mathcal{A} 
\end{align}
where $v_i$ is the velocity, $s_i$ the position of Agent $i$'s geometric center along the path and $u_i=a_i$ its control input constituting the longitudinal acceleration. By \prettyref{fig:problem_description_sketch}, the merging point (MP) is the origin of the curvilinear reference frame along $\mathcal{P}_i$, which is defined as \mbox{$s_i=s_{\text{MP}}\defeq0$.} Although being considered as a point mass in \eqref{eq:problem_modeling_sysContTime}, each \mbox{Agent $i$} exhibits a length $L_i \in \mathbb{R}_{>0}$ and width $W_i \in \mathbb{R}_{>0}$, which influences the  minimum safety distance in \prettyref{sec:NMPC_CA}.

To be used in an MPC-based control scheme, it is convenient to discretize the system dynamics \prettyref{eq:problem_modeling_sysContTime}. To this end, we apply zero-order hold (ZOH) discretization with sample time \mbox{$T_s \in \mathbb{R}_{>0}$} and obtain the exact discrete-time representation 
\begin{align} \label{eq:problem_modeling_sysDiscTime}
	\underbrace{\begin{bmatrix}
		{s}_{i,k+1} \\
		{v}_{i,k+1}
	\end{bmatrix}}_{x_{i,k+1}} =
    \begin{bmatrix}
    	1   & T_s \\
    	0   & 1
    \end{bmatrix} 
	\underbrace{\begin{bmatrix}
	{s}_{i,k} \\
	{v}_{i,k}
\end{bmatrix}}_{x_{i,k}} +
    \begin{bmatrix}
	\nicefrac{T_s^2}{2} \\
    T_s
\end{bmatrix}  \underbrace{\vphantom{
      \begin{bmatrix}
	{s}_{i,k} \\
	{v}_{i,k}
      \end{bmatrix}} a_{i,k}}_{u_{i,k}} 
\end{align}
with 
$x_{i,k} = [s_{i,k} ~\,v_{i,k}]^{\T}$,
$u_{i,k}=a_{i,k}$ at time $t_k=k T_s$, \mbox{$k\in\mathbb{N}_{0}$}. The control input $u_{i,k} \in \mathbb{U}_i$ is bounded by the compact set 
$\mathbb{U}_i \defeq [\underline{u}_i,\, \overline{u}_i]$ with $\underline{u}_i \leq 0 \leq \overline{u}_i$. 
With the lumped state and input vector \mbox{$x_{k} = [x_{1,k}^{\T} ~\,x_{2,k}^{\T}]^{\T} \in \mathbb{R}^{n_x}$} and 
\mbox{$u_{k} = [u_{1,k} ~\,u_{2,k}]^{\T} \in \mathbb{R}^{n_u}$},
respectively, where $n_x=4$ and $n_u=2$, the lumped system dynamics can concisely be written as $x_{k+1} = f(x_{k}, u_{k})$.

In the next section, we continue to formalize \prettyref{prob:problem_description_problemDef} in a centralized nonlinear MPC framework, in which collision avoidance constraints couple the agents while agent-related objectives and constraints are stated independently.

\section{NONLINEAR MPC FORMULATION}
\label{sec:NMPC}

When applying MPC, at every time $t_k$ we solve a finite-horizon optimal control problem (OCP) over a prediction horizon of \mbox{$N \in \mathbb{N}$} time steps. After optimization, only the first control input is applied to the plant and optimization is repeatedly executed over a shifted horizon at time $t_{k+1}$ \cite{Rawlings2022a}.

Hereafter, by $x_{j\mid k}$ we refer to the prediction of variable $x$ at the future time step $k+j$ given information up to time $k$, and by $x_{j\mid k}^\star$ to its optimal value. Moreover, for $a\leq b$, $c\leq d$, we denote by $\mathbb{N}_{[a,b]}$, $\mathbb{R}_{[c,d]}$ and $\mathbb{R}_{(c,d]}$ the sets \mbox{$\{ k \in \mathbb{N}_0 \mid a \leq k \leq b \}$,} $\{ k \in \mathbb{R} \mid c \leq k \leq d \}$ and $\{ k \in \mathbb{R} \mid c < k \leq d \}$, respectively.

\subsection{Agent Objectives and Constraints}
\label{sec:NMPC_agentObjCons}

As control objectives, every agent $i \in \mathcal{A}$ shall track its given reference velocity $v_{i,\text{ref}} \in \mathbb{R}_{>0}$, while at the same time minimizing the control input (acceleration) magnitude for reasons of comfort and efficiency. These objectives are captured through a quadratic stage cost at time step $k+j$
\begin{align*}
\ell_j(x_{j\mid k},u_{j\mid k}) &\defeq (x_{j\mid k} - x_{\text{ref},j\mid k})^{\T} Q\, (x_{j\mid k}-x_{\text{ref},j\mid k}) \\&~+ (u_{j\mid k})^{\T} R\, u_{j\mid k} 
\end{align*}
with reference state $x_{\text{ref},j\mid k} = [0~ v_{1,\text{ref},j\mid k}~ 0~ v_{2,\text{ref},j\mid k}]^{\T}$ (states $s_1, \,s_2$ are not controlled), and  terminal cost
\begin{align*}
\ell_N(x_{N\mid k}) \defeq (x_{N\mid k}-x_{\text{ref},N\mid k})^{\T} Q_N\, (x_{N\mid k}-x_{\text{ref},N\mid k})
\end{align*}
where $Q, \,Q_N,\succeq 0$ (weights on $s_1$, $s_2$ are zero) and $R \succeq 0$ are positive semi-definite, diagonal weight matrices.

The velocity of every agent shall be bounded between a minimum velocity $\underline{v}=0$ (only driving in forward direction) and a maximum velocity $\overline{v} \in \mathbb{R}_{>0}$. To this end, we define the  
constraint functions $h_{\underline{v}}$, $h_{\overline{v}}$ for Agent $i \in \mathcal{A}$, that is,
\begin{subequations}
\begin{align}
 	h_{\underline{v}}(x_{i,j\mid k})  &\defeq ~~v_{i,j\mid k} \label{eq:NMPC_agentObjCons_vmin}\\
	h_{\overline{v}}(x_{i,j\mid k})  &\defeq -v_{i,j\mid k} + \overline{v} \label{eq:NMPC_agentObjCons_vmax} 
\end{align}
\end{subequations}
and impose the following constraints over the horizon, i.e.,
\begin{align}
h_{\underline{v}}(x_{i,j\mid k}) \geq 0 ~~\land~~ h_{\overline{v}}(x_{i,j\mid k}) \geq 0,~~\forall j \in \mathbb{N}_{[1,N]}. \label{eq:NMPC_agentObjCons_vCons}
\end{align}
Moreover, the control input (acceleration) shall be bounded to accommodate physical actuator limitations, i.e.,
\begin{align*}
	u_{i,j\mid k} \in \mathbb{U}_i, ~~ \forall j \in \mathbb{N}_{[0,N-1]}. 
\end{align*}

\subsection{Safe Distance Between Agents}
\label{sec:NMPC_CA}

\begin{figure}[t!]
		\centering
		\setlength\fwidth{0.295\textwidth}	
		\input{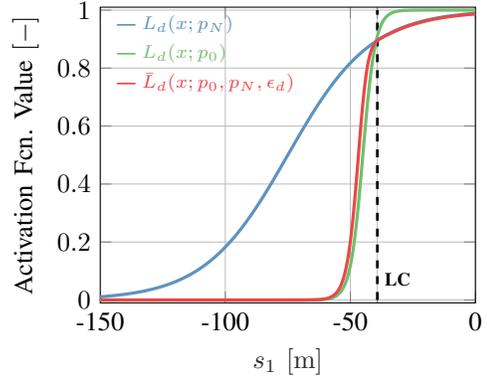}	
		\vspace*{-3mm}
		\caption{Activation function $L_d(x;p)$ to increase safety distance when Agent 1 approaches the lane change point, for $p_0=[0.4~\,{-45}]^{\T}$ (green) and $p_N=[0.06~\,{-75}]^{\T}$ (blue); $p_0$ enables overtaking much closer to LC. $\bar{L}_d(x;p_0,p_N)$ (red) is the activation function with $\epsilon_d = 0.0025$, which results from interpolation of $L_d(x;p_0)$ and $L_d(x;p_N)$. By construction, it holds $\bar{L}_d(x;p_0,p_N) \leq L_d(x;p_N)$ for all $s_1\in \mathbb{R}$, see \prettyref{sec:NMPCwCBF_safeDistHor}.} 
		\vspace*{-5mm}
		\label{fig:NMPC_CA_actFcnLd}
\end{figure}

To ease notation, we neglect time in the subsequent considerations. Referring to \prettyref{fig:problem_description_sketch}, the curvilinear reference frame of both agents has its origin at the merging point, i.e., at $s_1=s_2=s_{\text{MP}}\defeq 0$. That being said, we define the distance $d(x) \in \mathbb{R}_{\geqslant 0}$ between agents as the non-negative difference of path positions when projected onto a single path, i.e.,
\begin{align}
	d(x) \defeq \lvert s_1 - s_2 \rvert \geq d_s(x). 
	\label{eq:NMPC_CA_safeDistAbs}
\end{align}
To keep agents safe, $d(x)$ has to be lower bounded by a safety distance $d_{s}(x) \in \mathbb{R}_{\geqslant 0}$, which is sufficiently positive when Agent 1 moves into Agent 2's lane ($s_1 \geq s_{\text{LC}}$) and zero when both agents drive in separate lanes ($s_1 < s_{\text{LC}}$). That way, we avoid agents being forced to keep a safety distance before the lane change point is reached by \mbox{Agent 1}. Using the indicator function $\mathbf{1}_d(x)$, we can formalize \mbox{this condition as}
\begin{align*}
	d_s(x) \defeq 	\mathbf{1}_d(x) \,d_{\text{safe}}(x), ~~~\mathbf{1}_d(x) \defeq \begin{cases}
		1, & s_1 \geq s_{\text{LC}} \\
		0, & s_1 < s_{\text{LC}}
	\end{cases}
\end{align*} 
where $d_{\text{safe}}(x) \in \mathbb{R}_{> 0}$ is a state-dependent minimum safety distance. 
We define the minimum safety distance according to well-known practice in literature \cite{Li2011a} as
\begin{align} \label{eq:NMPC_CA_dsafe}
	d_{\text{safe}}(x) &\defeq d_0 + v_{f}(x) t_h
\end{align}
with stopping distance $d_0 \in \mathbb{R}_{>0}$ and time headway distance $v_{f}(x) t_h$, which depends on 
the headway time $t_h \in \mathbb{R}_{>0}$ and the velocity $v_f(x)$ of the following agent, i.e., 
\begin{align}	
	v_{f}(x) &\defeq \mathbf{1}_{lf}(x) \,v_1 + (1-\mathbf{1}_{lf}(x)) \,v_2. \label{eq:NMPC_CA_vf}
\end{align}
Here, $\mathbf{1}_{lf}(x)$ denotes the (\textit{leader/follower}) indicator function to determine, which agent is in front, i.e.,
\begin{align*}
	\mathbf{1}_{lf}(x) \defeq \begin{cases}
		1, & s_2 \geq s_1 \\
		0, & s_2 < s_1.
	\end{cases}
\end{align*} 
To ease computational complexity we aim to avoid mixed-integer optimization. To this end, we smoothly approximate  indicator functions $\mathbf{1}_{d}(x)$ and $\mathbf{1}_{lf}(x)$ by activation functions 
\begin{align}
	L_d(x;p) &\defeq \frac{1}{1+e^{-m_d(s_1-c_d)}}, ~~~p = \begin{bmatrix}
		m_d& c_d
	\end{bmatrix}^{\T} \label{eq:NMPC_CA_actFcnLd} \\
	\hspace*{-1mm}L_{lf}(x) &\defeq \frac{1}{1+e^{-m_{lf}(s_2-s_1)}} \label{eq:NMPC_CA_actFcnLlf}
\end{align}
with slopes $m_{lf},m_d \in \mathbb{R}_{>0}$ and offset $c_d\in\mathbb{R}$, see \prettyref{fig:NMPC_CA_actFcnLd} for two examples of $L_d(x;p)$. 
Starting from \prettyref{eq:NMPC_CA_safeDistAbs}, we define the constraint function $h_d$, that is,
\begin{align} \label{eq:NMPC_CA_CBFdist}
	h_d(x;p) \defeq (s_{1}-s_{2})^2 - ( \underbrace{L_d(x;p) \, \tilde{d}_{\text{safe}}(x)}_{\tilde{d}_s(x;p)} )^2,
\end{align}
where the smooth function $\tilde{d}_{\text{safe}}$ results from replacing $\mathbf{1}_{lf}(x)$ in \eqref{eq:NMPC_CA_vf} with $L_{lf}(x)$ in \eqref{eq:NMPC_CA_actFcnLlf}. 
We finally impose 
\begin{align*}
	h_d(x_{j\mid k};p) \geq 0, ~~\forall j \in \mathbb{N}_{[1,N]} 
\end{align*}
to ensure a safe distance between agents. If agents change their order before the lane change point (when projected onto a single path), we refer to this maneuver as \textit{overtaking}. 

\begin{remark} \normalfont
The parameter $p$ in $L_d(x;p)$ is chosen such that agents have at least $90\%$ of the minimum safety distance when Agent 1 reaches the lane change point, \mbox{see \prettyref{fig:NMPC_CA_actFcnLd}.}
\end{remark}

\subsection{Centralized Optimal Control Problem}
\label{sec:NMPC_COCP}

Given the state $x_k$ at time $t_k$, the resulting nonconvex OCP for lane merging can be stated as:

\vspace*{1mm}
\noindent
\textbf{Centralized Lane Merging NMPC Problem:}\\[-5mm]
\begin{subequations} \label{eq:NMPC_COCP_OCPdef}
	\begin{align}
		\underset{ U_k }{\mathmin} ~& \ell_N(x_{N\mid k})  +  \sum_{j=0}^{N-1} \ell_j(x_{j\mid k},u_{j\mid k}) \\
        \mathst~&~x_{j+1\mid k} = f(x_{j\mid k}, u_{j\mid k}),~~~~~~~~~~\forall j \in \mathbb{N}_{[0,N-1]}\\
        ~& ~h_{\underline{v}}(x_{i,j\mid k}) \geq 0,~~~~~~~~~ \forall i\in \mathcal{A} ,~\forall j \in \mathbb{N}_{[1,N]}\\
		~&~h_{\overline{v}}(x_{i,j\mid k}) \geq 0,~~~~~~~~~ \forall i\in \mathcal{A} ,~\forall j \in \mathbb{N}_{[1,N]}\\
		~&~h_{d}(x_{j\mid k}; p_0) \geq 0 ,~~~~~~~~~~~~~~~~~\forall j \in \mathbb{N}_{[1,N]}\\		
		~&~ u_{j\mid k} \in \mathbb{U}
        ,~~~~~~~~~~~~~~~~~~~~~~~~~~\forall j \in \mathbb{N}_{[0,N-1]}\\	
        ~&~x_{0\mid k} = x_k \\[-6mm]\notag
	\end{align}
\end{subequations}
where  $U_k \defeq (u_{0\mid k},u_{1\mid k}, \ldots, u_{N-1\mid k})$ is the control sequence to be optimized over the prediction horizon of length $N \in \mathbb{N}$, $p_0$ is a suitable parameter vector and $\mathbb{U}\defeq \mathbb{U}_1 \times \mathbb{U}_2$. 
In its current form, though, OCP \eqref{eq:NMPC_COCP_OCPdef} is not guaranteed to be recursively feasible  --- i.e., if OCP \eqref{eq:NMPC_COCP_OCPdef} is feasible at time $t_k$ there is no guarantee that it is feasible at time $t_{k+1}$.
Subsequently, we convey how to define terminal ingredients using (quasi-)DTCBF certificates to guarantee recursive feasibility.

\section{METHODOLOGY TO GUARANTEE RECURSIVE FEASIBILITY} 
\label{sec:DTCBF}

\subsection{Nonlinear MPC Setting}
\label{sec:DTCBF_background}

We consider discrete-time systems of the form \\[-5mm]
\begin{align} \label{eq:DTCBF_background_sysDiscTime}
	x_{k+1} = f(x_{k}, u_{k}) \\[-6mm]\notag
\end{align}
with state vector $x_k \in 
\mathbb{R}^{n_x}$, input vector \mbox{$u_k \in \mathbb{U} \subseteq \mathbb{R}^{n_u}$},  
vector field \mbox{$f: \mathbb{R}^{n_x} \times \mathbb{R}^{n_u} \rightarrow \mathbb{R}^{n_x}$}, where $n_x,n_u \in \mathbb{N}$, and input constraint set $\mathbb{U}$. 
\begin{assumption}[Model Uncertainties, Disturbances] \normalfont \label{ass:DTCBF_background_modelUncertainties}
The dynamical system \eqref{eq:DTCBF_background_sysDiscTime} is not affected by uncertainties or unmodeled exogenous disturbances. 
\end{assumption}
Without loss of generality, we consider the following nonlinear MPC problem:\\[-5mm]
\begin{subequations} \label{eq:DTCBF_background_NMPCsetup}
	\begin{align}
		\underset{ U_k }{\mathmin} ~& J(U_k; x_k) \label{eq:DTCBF_background_NMPCsetup_cost}\\
        \mathst~&~x_{j+1\mid k} = f(x_{j\mid k}, u_{j\mid k}),~~~~~~~~~~\forall j \in \mathbb{N}_{[0,N-1]} \label{eq:DTCBF_background_NMPCsetup_dynamics}\\
        ~& ~H(x_{j\mid k}) \geq 0,~~~~\,~~~~~~~~~~~~~~~~~\forall j \in \mathbb{N}_{[1,N-2]} \label{eq:DTCBF_background_NMPCsetup_constH}\\
        ~& ~h(x_{N-1\mid k}) \geq 0\label{eq:DTCBF_background_NMPCsetup_hNminus1}\\
        ~&~ h(x_{N\mid k}) \geq (1-\gamma) h(x_{N-1\mid k}) \label{eq:DTCBF_background_NMPCsetup_terminalh}\\
		~&~ u_{j\mid k} \in \mathbb{U}
        ,~~~~~~~~~~~~~~~~~~~~~~~~~~\forall j \in \mathbb{N}_{[0,N-1]}\label{eq:DTCBF_background_NMPCsetup_constU}\\	
        ~&~x_{0\mid k} = x_k \label{eq:DTCBF_background_NMPCsetup_initialCond} \\[-6mm]\notag
	\end{align}
\end{subequations}
where $J$ is a positive semi-definite cost, $H: \mathbb{R}^{n_x} \rightarrow \mathbb{R}$ a constraint function, $U_k \defeq (u_{0\mid k}, u_{1\mid k}, \ldots, u_{N-1\mid k})$ the control sequence to be optimized over the prediction horizon of length $N \in \mathbb{N}$ and $x_k$ the state at time $t_k$.

In line with \prettyref{sec:NMPC}, the objective is to find terminal ingredients \eqref{eq:DTCBF_background_NMPCsetup_hNminus1} and \eqref{eq:DTCBF_background_NMPCsetup_terminalh}, and thus an MPC-based control law that is recursively feasible and guarantees for all feasible initial states $x_0$ the satisfaction of state constraints \eqref{eq:DTCBF_background_NMPCsetup_constH}-\eqref{eq:DTCBF_background_NMPCsetup_terminalh} and input constraints \eqref{eq:DTCBF_background_NMPCsetup_constU}.
For fixed $\gamma \in \mathbb{R}_{(0,1]}$, the terminal ingredients are defined by the constraint function \mbox{$h: \mathbb{R}^{n_x} \rightarrow \mathbb{R}$} with 0-superlevel set \\[-5mm]
\begin{align} \label{eq:DTCBF_background_setC}
	\mathcal{C} \defeq \{ x \in \mathbb{R}^{n_x} \mid h(x) \geq 0  \}. \\[-6mm]\notag
\end{align}
\begin{theorem} \normalfont
\label{thm:DTCBF_background_consHorizon}
Suppose MPC \eqref{eq:DTCBF_background_NMPCsetup} is recursively feasible with $h$ and $H$ in \eqref{eq:DTCBF_background_NMPCsetup_constH}-\eqref{eq:DTCBF_background_NMPCsetup_terminalh} being equal.
Then, for same $h$, MPC \eqref{eq:DTCBF_background_NMPCsetup} is recursively feasible if $H$ is chosen (differently from $h$) such that it holds \mbox{$H(x) \geq h(x) \geq 0,~\forall x \in \mathcal{C}$}, or equivalently $\mathcal{C} \subseteq \mathcal{C}^\prime \defeq \{ x \in \mathbb{R}^{n_x} \mid H(x) \geq 0  \}$.
\end{theorem}
\vspace*{-2mm}
\begin{proof}
With $H$ being equal to $h$, terminal ingredients \eqref{eq:DTCBF_background_NMPCsetup_hNminus1}-\eqref{eq:DTCBF_background_NMPCsetup_terminalh} guarantee recursive feasibility of MPC \eqref{eq:DTCBF_background_NMPCsetup}. 
When choosing $H$ (differently from $h$), for recursive feasibility it has to hold: $h(x_{N-1\mid k}^\star)\geq 0$ implies $H(x_{N-2\mid k+1})\geq 0$. When setting $x_{N-2\mid k+1}=x_{N-1\mid k}^\star$, $h(x_{N-1\mid k}^\star)\geq 0$ and $H(x) \geq h(x) \geq 0$, \mbox{$\forall x \in \mathcal{C}$} imply $H(x_{N-2\mid k+1})\geq 0$. 
\hfill $\square$
%
%
%
\end{proof}
\vspace*{1mm}

\noindent
To determine $h$ in MPC \eqref{eq:DTCBF_background_NMPCsetup}, there are generally two options: 
\begin{enumerate}
    \item Define $h$ such that \eqref{eq:DTCBF_background_NMPCsetup} is recursively feasible with $H$ equal to $h$. Then, modify $H$ by respecting \prettyref{thm:DTCBF_background_consHorizon} to reduce conservatism (i.e., relax $\mathcal{C}=\mathcal{C}^\prime$ to $\mathcal{C} \subseteq \mathcal{C}^\prime$).
    \item Fix $H$, then define $h$ such that \eqref{eq:DTCBF_background_NMPCsetup} is recursively feasible (synthesis of $h$ needs to consider $H$). 
\end{enumerate}
In the remainder, we focus on the first option.
\begin{remark} \normalfont
By choosing $H$ different from $h$, the initial state can be located outside the safe set, i.e., \mbox{$x_{0\mid k} \notin \mathcal{C}$} 
and safely transition from $\mathcal{C}^\prime$ into $\mathcal{C}$ over the horizon. Thus, we reduce conservatism 
compared to \cite{Zeng2021a,Zeng2021b,Schillinger2021a}.
\end{remark}
\vspace*{-1mm}
\noindent
Evidently, \eqref{eq:DTCBF_background_NMPCsetup_terminalh} takes the form of a so called (quasi-)DTCBF certificate which requires certain conditions on $h$ to be satisfied, as further discussed in \prettyref{sec:DTCBF_terminalSetDTCBF} and \prettyref{sec:DTCBF_terminalSetTwoStep}.

\subsection{Terminal Ingredients: DTCBF Certificates}
\label{sec:DTCBF_terminalSetDTCBF}


\begin{definition}[Controlled invariance \cite{Ames2019a}] \normalfont
	\label{def:DTCBF_background_fwdInv}
The set $\mathcal{C}$ is said to be controlled invariant for the dynamical system \eqref{eq:DTCBF_background_sysDiscTime} with input constraint set $\mathbb{U}$, if for all $x_k \in \mathcal{C}$ there exists an input $u_k \in \mathbb{U}$ 
such that $x_{k+1} = f(x_k,u_k) \in \mathcal{C}$.
\end{definition}

\begin{definition}[Discrete-time CBF \cite{Agrawal2017a,Zeng2021a}] \normalfont
		\label{def:DTCBF_background_DTCBF}
	Let $\mathcal{C} \subset \mathbb{R}^{n_x}$ 
    be the 0-superlevel set of constraint function 
    \mbox{$h: \mathbb{R}^{n_x} \rightarrow \mathbb{R}$} 
    according to \prettyref{eq:DTCBF_background_setC}. For the dynamical system \prettyref{eq:DTCBF_background_sysDiscTime}, $h$ is a discrete-time control barrier function (DTCBF), 
    if there exists a \mbox{class-$\mathcal{K}_\infty$} function\footnote{A continuous function $\alpha: \mathbb{R}_{\geqslant 0} 
    \rightarrow \mathbb{R}_{\geqslant 0}$ 
    is said to belong to class-$\mathcal{K}_\infty$ if it is strictly increasing, $\alpha(0)=0$ and $\lim_{r\rightarrow \infty} \alpha(r)=\infty$.} \mbox{$\alpha: \mathbb{R}_{\geqslant 0} \rightarrow \mathbb{R}_{\geqslant 0}$} with $\alpha(r) \leq r$ for all \mbox{$r \in \mathbb{R}_{\geqslant 0}$} and if for all $x_k \in \mathcal{C}$ there exists a control input $u_k \in \mathbb{U}$ such that
	\vspace*{-1mm}
	\begin{align} \label{eq:DTCBF_background_DTCBFalpha}
            h(x_{k+1}) - h(x_k) \geq -\alpha(h(x_k)). \\[-6mm]\notag
	\end{align}
\end{definition}
\vspace*{1mm}

Inequality \eqref{eq:DTCBF_background_DTCBFalpha} with DTCBF $h$ is oftentimes referred to as DTCBF constraint \cite{Zeng2021a} or DTCBF certificate \cite{Taylor2022a}. Subsequently, we follow a frequently pursued approach in literature and define $\alpha$ as a linear function \cite{Agrawal2017a}, i.e., 
\vspace*{-1mm}
\begin{align*} 		
	h(x_{k+1}) \geq (1-\gamma) h(x_k), ~~ \gamma \in \mathbb{R}_{(0,1]}. 
\end{align*}

\begin{remark} \normalfont
Note that if $h$ is a DTCBF, then $\mathcal{C}$ is controlled invariant with input constraint set $\mathbb{U}$.
For $\gamma = 1$, we obtain controlled invariance in a \textit{classical} sense, while \mbox{$\gamma < 1$} constitutes a stricter controlled invariance condition.
\end{remark}

\begin{theorem}[Recursive Feasibility] \normalfont
\label{thm:DTCBF_terminalDTCBF_recFeasibility}
Suppose \prettyref{ass:DTCBF_background_modelUncertainties} holds, $h$ in \eqref{eq:DTCBF_background_NMPCsetup_hNminus1}  and \eqref{eq:DTCBF_background_NMPCsetup_terminalh} is a DTCBF and \mbox{$H(x)\geq h(x) \geq 0$,} $\forall x \in \mathcal{C}$ (cf. \prettyref{thm:DTCBF_background_consHorizon}). Then, MPC \eqref{eq:DTCBF_background_NMPCsetup} is recursively feasible, that is, for all feasible initial states $x_0$, MPC \eqref{eq:DTCBF_background_NMPCsetup} is guaranteed to remain feasible.
\end{theorem}
\begin{proof}
Suppose MPC \eqref{eq:DTCBF_background_NMPCsetup} with terminal ingredients \eqref{eq:DTCBF_background_NMPCsetup_hNminus1}  and \eqref{eq:DTCBF_background_NMPCsetup_terminalh} is initially feasible at time $t_0$. Then, at any time 
$t_{k} > t_0$, $k\in \mathbb{N}$ by \prettyref{def:DTCBF_background_DTCBF} $x_{N-1\mid k} = x_{N\mid k-1}^\star \in \mathcal{C}$ implies feasibility of \eqref{eq:DTCBF_background_NMPCsetup_terminalh} 
as $h$ is a DTCBF. That being said, MPC \eqref{eq:DTCBF_background_NMPCsetup} is recursively feasible. \hfill $\square$\\[-3.5mm] 
\end{proof} 


\subsection{Terminal Ingredients: quasi-DTCBF Certificates}
\label{sec:DTCBF_terminalSetTwoStep}

Subsequently, we propose an alternative method to define terminal ingredients \eqref{eq:DTCBF_background_NMPCsetup_hNminus1} and \eqref{eq:DTCBF_background_NMPCsetup_terminalh} for recursive feasibility. \\[-5mm]
\begin{definition}[quasi-DTCBF, Two-step Condition] \normalfont
		\label{def:DTCBF_terminalSetTwoStep_twoStepCond}
Let \mbox{$\mathcal{C} \subset \mathbb{R}^{n_x}$} be the \mbox{0-superlevel} set of constraint function \mbox{$h: \mathbb{R}^{n_x} \rightarrow \mathbb{R}$} as in \prettyref{eq:DTCBF_background_setC}. For the dynamical system \prettyref{eq:DTCBF_background_sysDiscTime}, $h$ is a quasi-DTCBF (qDTCBF) if for given $\gamma \in \mathbb{R}_{(0,1]}$ and for all $x_k \in \mathcal{C}$ and all $u_k \in \mathbb{U}$ with 
\\[-5mm]
\begin{align}
\label{eq:DTCBF_terminalSetTwoStep_DTCFBlikeCons}
h(f(x_k,u_k)) \geq (1-\gamma) h(x_{k}) \\[-6mm]\notag
\end{align}
there exists $u_{k+1} \in \mathbb{U}$ such that \\[-5mm]
\begin{align*}
h(f(x_{k+1},u_{k+1}) \geq (1-\gamma) h(x_{k+1}) \\[-6mm]
\end{align*}
where $x_{k+1}=f(x_k,u_k)$ (referred to as \textit{two-step condition}). 
\end{definition}
\begin{lemma} \normalfont
\label{lem:DTCBF_terminalSetTwoStep_CzeroCI}
Following  \prettyref{def:DTCBF_terminalSetTwoStep_twoStepCond}, we define the set \mbox{$\mathcal{C}_0 \defeq \{ x \in \mathcal{C} \mid \exists u \in \mathbb{U} : h(f(x,u)) \geq (1-\gamma) h(x) \}$.} If $h$ is a qDTCBF, $\mathcal{C}_0$ is controlled invariant.
\end{lemma}
\vspace*{-1mm}
\begin{proof}
Let $x_k \in \mathcal{C}_0$, then $h(x_k) \geq 0$ and there is $u_k \in \mathbb{U}$ such that $h(f(x_k,u_k)) \geq (1-\gamma) h(x_k)$ and thus $x_{k+1}=f(x_k,u_k) \in \mathcal{C}$. Since $h$ is a qDTCBF, for $x_{k+1} \in \mathcal{C}$ there is $u_{k+1} \in \mathbb{U}$ such that $h(f(x_{k+1},u_{k+1}) \geq (1-\gamma) h(x_{k+1})$. Hence, $x_{k+1} \in \mathcal{C}_0$ and $\mathcal{C}_0$ is controlled invariant. \hfill $\square$
\end{proof}
\vspace*{1mm}

Hereafter, we refer to \eqref{eq:DTCBF_terminalSetTwoStep_DTCFBlikeCons} as qDTCBF certificate. 
Compared to \prettyref{def:DTCBF_background_DTCBF}, 
\prettyref{def:DTCBF_terminalSetTwoStep_twoStepCond} is a weaker condition 
as it only requires $h(x_{k+1}) \geq (1-\gamma) h(x_{k})$ to be feasible for all $x_{k}\in \mathcal{C}_0 \subseteq \mathcal{C}$. Thus, $\mathcal{C}$ is not necessarily controlled invariant unless $\mathcal{C}_0=\mathcal{C}$. 
If $h$ is a DTCBF, though, then $h$ also satisfies \prettyref{def:DTCBF_terminalSetTwoStep_twoStepCond} as $\mathcal{C}_0=\mathcal{C}$ and $x_k \in \mathcal{C}$ implies feasibility of \eqref{eq:DTCBF_terminalSetTwoStep_DTCFBlikeCons} by \prettyref{def:DTCBF_background_DTCBF} and thus feasibility of \eqref{eq:DTCBF_terminalSetTwoStep_DTCFBlikeCons} at the next time instance. The converse, however, does not hold.
\begin{theorem}[Recursive Feasibility] \normalfont
		\label{thm:DTCBF_terminalSetTwoStep_recFeasibility}
   Suppose \prettyref{ass:DTCBF_background_modelUncertainties} holds, $h$ in \eqref{eq:DTCBF_background_NMPCsetup_hNminus1}  and \eqref{eq:DTCBF_background_NMPCsetup_terminalh} is a qDTCBF and \mbox{$H(x)\geq h(x) \geq 0$,} $\forall x \in \mathcal{C}$ (cf. \prettyref{thm:DTCBF_background_consHorizon}). 
   Then, MPC \eqref{eq:DTCBF_background_NMPCsetup} is recursively feasible, that is, for all feasible initial states $x_0$, MPC \eqref{eq:DTCBF_background_NMPCsetup} is guaranteed to remain feasible. 
\end{theorem}
\vspace*{-1mm}
\begin{proof}
Suppose MPC \eqref{eq:DTCBF_background_NMPCsetup} with terminal ingredients \eqref{eq:DTCBF_background_NMPCsetup_hNminus1}  and \eqref{eq:DTCBF_background_NMPCsetup_terminalh} is initially feasible at time $t_0$ such that by \prettyref{def:DTCBF_terminalSetTwoStep_twoStepCond} it holds $x_{N-1\mid 0}^\star \in \mathcal{C}_0$. Then, at any time $t_{k} > t_0$, $k\in \mathbb{N}$ by \prettyref{def:DTCBF_terminalSetTwoStep_twoStepCond} and 
\mbox{\prettyref{lem:DTCBF_terminalSetTwoStep_CzeroCI}} $x_{N-1\mid k-1}^\star \in \mathcal{C}_0$ and feasibility of \eqref{eq:DTCBF_background_NMPCsetup_terminalh} at time $t_{k-1}$ implies feasibility of \eqref{eq:DTCBF_background_NMPCsetup_terminalh} at time $t_k$. 
As such, MPC \eqref{eq:DTCBF_background_NMPCsetup} is recursively feasible. \hfill $\square$
\end{proof}

\subsection{Properties of the Terminal Ingredients}
\label{sec:DTCBF_commonProps}


The advantage of qDTCBFs is that \prettyref{def:DTCBF_terminalSetTwoStep_twoStepCond} can be satisfied with lower absolute input bounds $\lvert \overline{u} \rvert$, $\lvert \underline{u} \rvert$ (element-wise if $\dim(u)>1$) than \prettyref{def:DTCBF_background_DTCBF} (DTCBFs) if \mbox{$\mathcal{C}_0 \subset \mathcal{C}$}, see 
\prettyref{sec:results_costReduction}. This is highly advantageous if only limited actuation is possible. For recursive feasibility, though, \prettyref{thm:DTCBF_terminalSetTwoStep_recFeasibility} requires $x_{N-1\mid k}$ to reside in $\mathcal{C}_0 \subseteq \mathcal{C}$ while \prettyref{thm:DTCBF_terminalDTCBF_recFeasibility} allows $x_{N-1\mid k}$ to be located in the larger \mbox{set $\mathcal{C}$.} 

When considering \mbox{$\mathbb{X}_{N}(x_{N-1\mid k}) \defeq \{ \,x_{N\mid k} \in \mathbb{R}^{n_x} \mid \eqref{eq:DTCBF_background_NMPCsetup_terminalh}\, \}$,} 
for both, DTCBF and qDTCBF certificates, $\mathbb{X}_{N}(x_{N-1\mid k})$ is a contracted subset of $\mathcal{C}$ (unless $h(x_{N-1\mid k})=0$ or $\gamma=1$, then $\mathbb{X}_{N}=\mathcal{C}$),
where the contraction is driven by $\gamma \in \mathbb{R}_{(0,1]}$ (lower $\gamma$ means stronger contraction). This additional tuning knob encourages an earlier reaction of the control system and reduces the cumulative MPC 
cost, see \prettyref{sec:results_costReduction}. At the same time, though, the contraction leads to a stricter terminal ingredient \eqref{eq:DTCBF_background_NMPCsetup_terminalh}.
Hence, the choice of $\gamma$ has to ensure 
a suitable \mbox{trade-off between contraction and cost \cite{Zeng2021a,Zeng2021b}.} 

If $h$ in \eqref{eq:DTCBF_background_NMPCsetup_hNminus1} and \eqref{eq:DTCBF_background_NMPCsetup_terminalh} is a vector-valued function such that $h(x) = (h_\iota(x))_{\iota\in\mathbb{N}_{[1,M]}}$, $M \in \mathbb{N}$ with  
$h_\iota : \mathbb{R}^{n_x} \rightarrow \mathbb{R}$ 
(in \eqref{eq:DTCBF_background_NMPCsetup} we assume \mbox{$M=1$),} joint feasibility of multiple (q)DTCBF certificates 
subject to input constraints has additionally to be guaranteed. This part of the theory is not in the scope of this paper. \prettyref{sec:NMPCwCBF}, though, conveys how joint feasibility can be achieved for the given lane merging application. 
\section{LANE MERGING NMPC WITH RECURSIVE FEASIBILITY GUARANTEES}
\label{sec:NMPCwCBF}

The methodology in \prettyref{sec:DTCBF} can be applied to discrete-time systems of the form \eqref{eq:DTCBF_background_sysDiscTime}. Nonlinear sampled-data systems can still be handled through an additional design step 
\cite{Schillinger2021a,Taylor2022a}.
In our application, the system \eqref{eq:problem_modeling_sysContTime} is linear and we obtain an exact discretization of the form \eqref{eq:DTCBF_background_sysDiscTime} using ZOH. 
To form terminal ingredients, we design $h$ first \mbox{(option 1} in \prettyref{sec:DTCBF_background}) and then adapt $H$ to reduce conservatism.

\subsection{Terminal Velocity DTCBF Certificates}
\label{sec:NMPCwCBF_agentCons}
By virtue of \prettyref{sec:NMPC_agentObjCons}, 
every agent needs to satisfy velocity constraints \eqref{eq:NMPC_agentObjCons_vCons} subject to input constraints \mbox{$u_{i} \in \mathbb{U}_i$}. We construct DTCBF certificates from 
\mbox{\eqref{eq:NMPC_agentObjCons_vmin}-\eqref{eq:NMPC_agentObjCons_vmax}}, i.e.,
\begin{subequations}
\begin{align}
	h_{\underline{v}}(x_{i,k+1}) &\geq (1-\gamma_v) h_{\underline{v}}(x_{i,k}) \label{eq:NMPCwCBF_agentCons_vMinCBFCons},  \\ 
	h_{\overline{v}}(x_{i,k+1})   &\geq (1-\gamma_v) h_{\overline{v}}(x_{i,k}),~\gamma_v \in \mathbb{R}_{(0,1]}.    \label{eq:NMPCwCBF_agentCons_vMaxCBFCons}
\end{align}
\end{subequations}

\begin{proposition} \normalfont
\label{prop:NMPCwCBF_agentCons_propVelCBF}
  Constraint functions \mbox{\eqref{eq:NMPC_agentObjCons_vmin}-\eqref{eq:NMPC_agentObjCons_vmax}} are DTCBFs on their respective 0-superlevel sets $\mathcal{C}_{\underline{v}}$, $\mathcal{C}_{\overline{v}}$ for the dynamical system \eqref{eq:problem_modeling_sysDiscTime} in accordance with \prettyref{def:DTCBF_background_DTCBF}. Moreover, for all $x_{i,k} \in \mathcal{C}_{\underline{v}} \cap \mathcal{C}_{\overline{v}}$ there exists a control input $u_{i,k} \in \mathbb{U}_i$ such that DTCBF certificates \mbox{\eqref{eq:NMPCwCBF_agentCons_vMinCBFCons}-\eqref{eq:NMPCwCBF_agentCons_vMaxCBFCons}} are jointly feasible. 
\end{proposition}
\vspace*{-1mm}
\begin{proof}
	Substituting dynamics \eqref{eq:problem_modeling_sysDiscTime} into \eqref{eq:NMPCwCBF_agentCons_vMinCBFCons}-\eqref{eq:NMPCwCBF_agentCons_vMaxCBFCons} yields
    \vspace*{-1mm}
	\begin{align*}
        {u}_{i,k}  \geq  -\frac{\gamma_v}{T_s} v_{i,k}  \,~~~~&\land~~~\,
        {u}_{i,k}  \leq \frac{\gamma_v}{T_s}  (\overline{v}-v_{i,k}), \notag\\[-7mm]
	\end{align*}
	 which translates into a non-positive lower bound and a non-negative upper bound on the input. 
	 Imposing these constraints and $u_{i,k} \in \mathbb{U}_i$ simultaneously is equivalent to
	 \begin{align} \label{eq:NMPCwCBF_agentCons_umaxmin}
	 	\hspace*{-2mm}\max\left\{ \underline{u},\, -\frac{\gamma_v}{T_s} v_{i,k}\right\}  \leq {u}_{i,k} \leq 
	 	\min\left\{ \overline{u},\, \frac{\gamma_v}{T_s}  (\overline{v}-v_{i,k}) \right\}
	 \end{align}
       for $0 \leq v_{i,k} \leq \overline{v}$, 
	   which implies joint feasibiliy of \mbox{\eqref{eq:NMPCwCBF_agentCons_vMinCBFCons}-\eqref{eq:NMPCwCBF_agentCons_vMaxCBFCons}} subject to $u_{i,k} \in \mathbb{U}_i$. \hfill$\square$ 
\end{proof}
\begin{remark} \normalfont
By virtue of \eqref{eq:NMPCwCBF_agentCons_umaxmin}, terminal velocity DTCBF certificates do not introduce any noticeable conservatism. Hence, we do not compare them with qDTCBF certificates.
\end{remark}

\subsection{Terminal Safety Distance qDTCBF Certificate}
\label{sec:NMPCwCBF_safeDist}

Applying qDTCBF certificates to guarantee safety offers advantages over a DTCBF formulation as shown in the remainder. By virtue of \eqref{eq:NMPC_CA_CBFdist}, we form the qDTCBF certificate \\[-5mm]
\begin{align}  
	h_d(x_{k+1}; p_N) &\geq (1-\gamma_d) h_{d}(x_{k}; p_N),~~\gamma_d \in \mathbb{R}_{(0,1]} \label{eq:NMPCwCBF_safeDist_distCBFCons} \\[-6mm] \notag
\end{align}
with parameter vector \mbox{$p_N = [m_{d,N}~ c_{d,N}]^{\T}$}. To ensure compliance of $h_d$ with 
\prettyref{def:DTCBF_terminalSetTwoStep_twoStepCond}, we solve the following verification problem for given $p_N$, $\gamma_d$, input constraint set $\mathbb{U}$ and fixed feedback policy $u_k=\kappa(x_k)$:\\[1mm]
\noindent
\textbf{Verification Problem for qDTCBF $h_d$ (cf. \prettyref{def:DTCBF_terminalSetTwoStep_twoStepCond}):}\\[-5.5mm]
\begin{subequations} \label{eq:NMPCwCBF_safeDist_verifyTwoStep}
	\begin{align}
      \hspace*{-1.8mm}0 \leq \min_{x_k \in \mathbb{R}^{n_x}}  
      &~ h_d(x_{k+2}; p_N) - (1-\gamma_d) h_{d}(x_{k+1}; p_N)  \label{eq:NMPCwCBF_safeDist_verifyTwoStep_step2} \raisetag{6mm}\\
      \hspace*{-1.8mm}\mathst~ &~h_d(x_{k+1}; p_N) - (1-\gamma_d) h_{d}(x_{k}; p_N) \geq 0 \label{eq:NMPCwCBF_safeDist_verifyTwoStep_step1}\\
      &~h_d(x_k;p_N) \geq 0 \label{eq:NMPCwCBF_safeDist_verifyTwoStep_xinC}\\
      &~\Delta v_k \geq \Delta \underline{v} \label{eq:NMPCwCBF_safeDist_verifyTwoStep_deltaV}
    \end{align}
\end{subequations}
with $x_{k+1}=f(x_k,\kappa(x_k))$ and $x_{k+2}=f(x_{k+1},\kappa(x_{k+1}))$, claiming that $x_k \in \mathcal{C}_0$ (cf. \eqref{eq:NMPCwCBF_safeDist_verifyTwoStep_step1}-\eqref{eq:NMPCwCBF_safeDist_verifyTwoStep_xinC}) implies feasibility of \eqref{eq:NMPCwCBF_safeDist_verifyTwoStep_step2}. The additional constraint \eqref{eq:NMPCwCBF_safeDist_verifyTwoStep_deltaV} on the relative velocity 
\begin{align*}
\Delta v_{k} \defeq L_{lf}(x_{k})( v_{2,k} - v_{1,k}) + (1-L_{lf}(x_{k}))( v_{1,k} - v_{2,k} ) 
\end{align*}
with $\Delta \underline{v} \in \mathbb{R}$, depending on the order of agents, simplifies to obtain a feasible solution of \eqref{eq:NMPCwCBF_safeDist_verifyTwoStep}. As feedback policy (to form inputs $u_k, u_{k+1}$ 
to satisfy \prettyref{def:DTCBF_terminalSetTwoStep_twoStepCond}), we choose
\vspace*{-1mm}
\begin{subequations} \label{eq:NMPCwCBF_safeDist_policy}
\begin{align} 
\kappa(x_k) &\defeq [ \tilde{a}_{1,k} ~ \tilde{a}_{2,k}]^{\T}~~\text{with} \label{eq:NMPCwCBF_safeDist_policy_kx}\\
\tilde{a}_{1,k} &\defeq  L_{lf}(x_{k}) \max \bigl\{ \underline{u},\, -\frac{\gamma_v}{T_s} v_{1,k}\bigr\} ~+ \label{eq:NMPCwCBF_safeDist_policy_input1}\\
&~~~~(1-L_{lf}(x_{k})) \min\bigl\{ \overline{u},\, \frac{\gamma_v}{T_s}  (\overline{v}-v_{1,k}) \bigr\} \notag\\
\tilde{a}_{2,k} &\defeq  L_{lf}(x_{k}) \min\bigl\{ \overline{u},\, \frac{\gamma_v}{T_s}  (\overline{v}-v_{2,k}) \bigr\}  ~+ \label{eq:NMPCwCBF_safeDist_policy_input2}\\
&~~~~(1-L_{lf}(x_{k})) \max \bigl\{ \underline{u},\, -\frac{\gamma_v}{T_s} v_{2,k}\bigr\} \notag
\end{align}
\end{subequations}
which enforces the leading agent to accelerate with maximum acceleration and the following agent to decelerate with maximum deceleration consistent with \eqref{eq:NMPCwCBF_agentCons_umaxmin} --- the most effective action to increase the distance between agents. 

We solve the nonconvex problem \eqref{eq:NMPCwCBF_safeDist_verifyTwoStep} with a global optimization algorithm to ensure that the minimal cost is globally greater than zero such as to avoid false conclusions. To this end, we utilize our variant of the $\alpha$BB algorithm \cite{Androulakis1995a}, which is not in the scope of this contribution. If \eqref{eq:NMPCwCBF_safeDist_verifyTwoStep} is feasible, a valid certificate is found, otherwise $\gamma_d$, $p_N$, $\Delta \underline{v}$ and/or $\mathbb{U}$ have to be modified until \eqref{eq:NMPCwCBF_safeDist_verifyTwoStep} is feasible. 
\begin{proposition} \normalfont
\label{prop:NMPCwCBF_safeDist_propSafeDistCBF}
Let $(p_N,\, \gamma_d,\, \Delta \underline{v})$ be parameters
for which verification problem \eqref{eq:NMPCwCBF_safeDist_verifyTwoStep} is feasible. Then, qDTCBF certificate \eqref{eq:NMPCwCBF_safeDist_distCBFCons} is jointly feasible with DTCBF certificates \mbox{\eqref{eq:NMPCwCBF_agentCons_vMinCBFCons}-\eqref{eq:NMPCwCBF_agentCons_vMaxCBFCons}} subject to input constraints $u \in \mathbb{U}$.
\end{proposition}
\vspace*{-1mm}
\begin{proof}
Joint feasibilty of \eqref{eq:NMPCwCBF_safeDist_distCBFCons} and \eqref{eq:NMPCwCBF_agentCons_vMinCBFCons}-\eqref{eq:NMPCwCBF_agentCons_vMaxCBFCons} subject to input constraints follows from the choice of the feedback policy \eqref{eq:NMPCwCBF_safeDist_policy} such as to satisfy \eqref{eq:NMPCwCBF_agentCons_umaxmin}.  \hfill $\square$
\end{proof}

\subsection{Safety Distance Constraints Over Remaining Horizon}
\label{sec:NMPCwCBF_safeDistHor}

Referring to MPC \eqref{eq:DTCBF_background_NMPCsetup} 
in \prettyref{sec:DTCBF_background}, we use the same velocity constraint functions $h_{\underline{v}}$ and $h_{\overline{v}}$ to form terminal ingredients and the constraints over the horizon as no conservatism is introduced that way. Conversely, imposing $h_d(x_{j\mid k};p_N)\geq 0$ with $p_N=[m_{d,N}~c_{d,N}]^{\T}$ for $j\in\mathbb{N}_{[1,N-2]}$ may cause significant conservatism as the slope $m_{d,N}$ of the activation function has to be chosen small to obtain a feasible solution of verification problem \eqref{eq:NMPCwCBF_safeDist_verifyTwoStep}. That way, the minimum safety distance between agents would be increased already far away from the lane change point, see $L_d(x;p_N)$ (blue) in \prettyref{fig:NMPC_CA_actFcnLd}, which may prevent agents from overtaking. 


Hence, for given $h_d$ with parameter $p_N=[m_{d,N}~c_{d,N}]^{\T}$ and 0-superlevel set $\mathcal{C}_d$ we aim to find a constraint function $H_d$ with parameter $p_0=[m_{d,0}~c_{d,0}]^{\T}$, which exhibits a higher slope \mbox{$m_{d,0} > m_{d,N}$} and satisfies \mbox{$H_d(x;p_0)\geq h_d(x;p_N) \geq 0$} for all $x \in \mathcal{C}_d$ (cf. \prettyref{thm:DTCBF_background_consHorizon}). Choosing $H_d$ equal to $h_d$ with parameter $p_0$ does not satisfy \prettyref{thm:DTCBF_background_consHorizon} as it requires $L_d(x;p_0)\leq L_d(x;p_N)$ to hold for all $s_1 \in \mathbb{R}$,
see $L_d(x;p_0)$ (green) in \prettyref{fig:NMPC_CA_actFcnLd}.   
Thus, we define the \textit{interpolated} 
activation function $\bar{L}_{d}$ (red curve in \prettyref{fig:NMPC_CA_actFcnLd}):\\[-5.5mm]
\begin{align*}
\bar{L}_{d}(x; p_0, p_N) \defeq L_d(x; p_0) \bigl(1 \hspace*{-0.3mm}+\hspace*{-0.3mm} L_d(x; p_N) \hspace*{-0.3mm}-\hspace*{-0.3mm} L_d(x; p_0) \hspace*{-0.3mm}-\hspace*{-0.3mm} \epsilon_d \bigr) \\[-5.5mm]
\end{align*}
with small $\epsilon_d \in \mathbb{R}_{>0}$ (in our case $\epsilon_d=0.0025$) to satisfy $\bar{L}_{d}(x; p_0, p_N)\leq L_d(x; p_N)$ for all $s_1 \in \mathbb{R}$,
and define the minimum safety distance \\[-5.5mm]
\begin{align}
\label{eq:NMPCwCBF_safeDist_safeDist}
\bar{d}_s(x;p_0,p_N) \defeq \bar{L}_{d}(x; p_0, p_N) \, \tilde{d}_{\text{safe}}(x) \\[-5.5mm]\notag
\end{align}
to devise the less conservative constraint function $H_d$, i.e., \\[-5.5mm]
\begin{align*}
    \hspace*{-2mm} H_d(x; p_0, p_N) &\defeq (s_{1}-s_{2})^2 - \bar{d}_s(x;p_0,p_N)^2 \\[-5.5mm]
\end{align*}
and impose $H_d(x_{j\mid k}; p_0, p_N) \geq 0$ for all $j\in\mathbb{N}_{[1,N-2]}$ to allow the agents to overtake close to the lane change point. 

\begin{figure*}[ht!]
		\centering
		\setlength\fwidth{0.84\textwidth}	
		\input{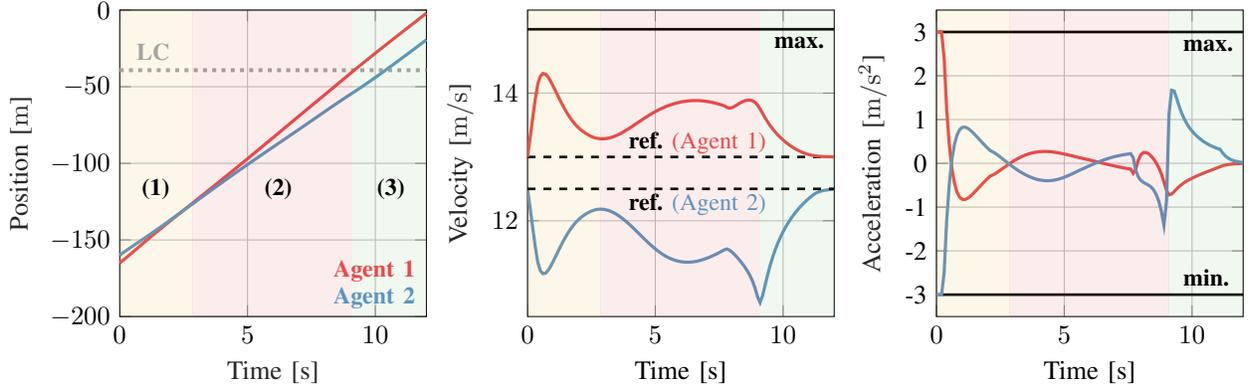}	
		\vspace*{-4mm}
		\caption{Validation Scenario: Agent 1 (red) overtakes Agent 2 (blue) before the lane change (LC) point. The maneuver can be subdivided into 3 phases: \mbox{(1) overtaking}, (2) accommodate increasing safe distance and (3) velocity reference tracking. State and input constraints are satisfied at all times.} 
		\vspace*{-5mm}
		\label{fig:sim_results_overtakingPosVelAccel}
\end{figure*}

\subsection{Resulting Recursively Feasible OCP}
\label{sec:NMPCwCBF_OCP}


With DTCBF certificates \eqref{eq:NMPCwCBF_agentCons_vMinCBFCons}, \eqref{eq:NMPCwCBF_agentCons_vMaxCBFCons} and qDTCBF certificate \eqref{eq:NMPCwCBF_safeDist_distCBFCons}, the (q)DTCBF-NMPC problem with recursive feasibility guarantees finally reads:
\vspace{1mm}

\noindent
\textbf{Centralized Lane Merging (q)DTCBF-NMPC Problem:} \\[-5mm]
\begin{subequations} \label{eq:NMPCwCBF_OCP_defOCP}
	\begin{align}
		\underset{ U_k }{\mathmin} ~&\, \ell_N(x_{N\mid k})  +  \sum_{j=0}^{N-1} \ell_j(x_{j\mid k},u_{j\mid k}) \label{eq:NMPCwCBF_OCP_defOCP_cost}\\
        \mathst~&\, x_{j+1\mid k} = f(x_{j\mid k}, u_{j\mid k}),~~~~~~~~~~\,\forall j \in \mathbb{N}_{[0,N-1]}\\  
  		~&\,H_{d}(x_{j\mid k}; p_0, p_N) \geq 0 ,~~~~~~\,\,~~~~~\forall j \in \mathbb{N}_{[1,N-2]} \label{eq:NMPCwCBF_OCP_defOCP_safeDistNminus2}\\
		~&\, h_{d}(x_{N-1\mid k}; p_N) \geq 0 \label{eq:NMPCwCBF_OCP_defOCP_safeDistNminus1} \\
        ~&\, h_{\underline{v}}(x_{i,j\mid k}) \geq 0,~~~~~~\,~~~ \forall i\in \mathcal{A} ,~\forall j \in \mathbb{N}_{[1,N-1]} \label{eq:NMPCwCBF_OCP_defOCP_vmin}\\
		~&\,h_{\overline{v}}(x_{i,j\mid k}) \geq 0,~~~~~\,~~~~ \forall i\in \mathcal{A} ,~\forall j \in \mathbb{N}_{[1,N-1]} \label{eq:NMPCwCBF_OCP_defOCP_vmax}\\
	    ~&\, h_{d}(x_{N\mid k}; p_N) \geq (1-\gamma_{d}) h_{d}(x_{N-1\mid k}; p_N) \label{eq:NMPCwCBF_OCP_defOCP_termSetSafety}\\
	    ~&\, h_{\underline{v}}(x_{i,N\mid k}) \,~~~\geq (1-\gamma_v) h_{\underline{v}}(x_{i,N-1\mid k}),\, \hspace*{-0.1mm}\forall i\in \mathcal{A} \label{eq:NMPCwCBF_OCP_defOCP_termSetVMin} \raisetag{4mm}\\
	    ~&\, h_{\overline{v}}(x_{i,N\mid k}) \,~~~\geq (1-\gamma_v) h_{\overline{v}}(x_{i,N-1\mid k}),\, \forall i\in \mathcal{A} \label{eq:NMPCwCBF_OCP_defOCP_termSetVMax} \\
        ~&\, \Delta v_{N-1\mid k} \geq \Delta \underline{v}  \label{eq:NMPCwCBF_OCP_defOCP_deltaV}\\	     
        ~&\, u_{j\mid k} \in \mathbb{U} ,~~~~~~~~~~~~~~~~~~\,~~~~~~~~\forall j \in \mathbb{N}_{[0,N-1]}\label{eq:NMPCwCBF_OCP_defOCP_inputCons} \\
        ~&\, x_{0\mid k} = x_{k} \\[-6mm]\notag
	\end{align}
\end{subequations}
with control input sequence $U_k \defeq (u_{0\mid k},u_{1\mid k}, \ldots, u_{N-1\mid k})$. 
Importantly, we need to impose \eqref{eq:NMPCwCBF_OCP_defOCP_deltaV} as this condition has also been applied in verification problem \eqref{eq:NMPCwCBF_safeDist_verifyTwoStep}. 
\begin{theorem} \normalfont
\label{thm:NMPCwCBF_OCP_recFeasible}
    Suppose \prettyref{ass:DTCBF_background_modelUncertainties} holds. Then, the centralized (q)DTCBF-NMPC \eqref{eq:NMPCwCBF_OCP_defOCP} is recursively feasible. 
\end{theorem}
\vspace*{-1mm}
\begin{proof}
    First, the definition of constraints \mbox{\eqref{eq:NMPCwCBF_OCP_defOCP_safeDistNminus2}-\eqref{eq:NMPCwCBF_OCP_defOCP_vmax}} follows \prettyref{thm:DTCBF_background_consHorizon}. Second, suppose $u_{0\mid k}^\star,\ldots,u_{N-1\mid k}^\star$ and $x_{1\mid k}^\star,\ldots,x_{N\mid k}^\star$ are the optimal input and state sequence at time $t_k$. We can always satisfy \eqref{eq:NMPCwCBF_OCP_defOCP_deltaV} 
    if \mbox{OCP \eqref{eq:NMPCwCBF_OCP_defOCP}} is feasible at time $t_k$: \textsc{i}) if $\Delta v_{N\mid k}^\star \geq \Delta \underline{v}$, we can choose $u_{N-2\mid k+1}=u_{N-1\mid k}^\star$; \textsc{ii}) if $\Delta v_{N\mid k}^\star < \Delta \underline{v}$, there exists $u_{N-2\mid k+1} \in \mathbb{U}$ (e.g., $u_{N-2\mid k+1}=-u_{N-1\mid k}^\star$ as \mbox{$\Delta v_{N-1\mid k}^\star \geq \Delta \underline{v}$)} to also satisfy \mbox{$h_{d}(x_{N-1\mid k+1}; p_N) \geq 0$} with even higher relative velocity. 
    Finally, 
    recursive feasibility of OCP \eqref{eq:NMPCwCBF_OCP_defOCP} follows from \prettyref{thm:DTCBF_terminalSetTwoStep_recFeasibility} together with \prettyref{prop:NMPCwCBF_agentCons_propVelCBF} and \prettyref{prop:NMPCwCBF_safeDist_propSafeDistCBF}, which include the proof of joint feasibility of 
    \eqref{eq:NMPCwCBF_OCP_defOCP_termSetSafety}-\eqref{eq:NMPCwCBF_OCP_defOCP_termSetVMax}
    subject to input \mbox{constraints \eqref{eq:NMPCwCBF_OCP_defOCP_inputCons}.} \hfill $\square$
\end{proof}

\section{SIMULATION RESULTS}
\label{sec:results}

\begin{figure}[b!]
			\centering
		\vspace*{-4mm}
		\setlength\fwidth{0.27\textwidth}	
%
%

\definecolor{mycolor1}{rgb}{0.92900,0.69400,0.12500}%
\definecolor{mycolor2}{RGB}{233,72,73}%
\definecolor{mycolor3}{RGB}{113,191,110}%

\begin{tikzpicture}

\begin{axis}[%
width=0.951\fwidth,
height=0.75\fwidth,
at={(0\fwidth,0\fwidth)},
scale only axis,
xmin=0,
xmax=12,
xlabel style={font=\color{white!15!black}},
xlabel={Time [$\mathrm{s}$]},
ymin=-0.4,
ymax=18,
ylabel style={font=\color{white!15!black}},
ylabel={Distance [$\mathrm{m}$]},
axis background/.style={fill=white},
xmajorgrids,
ymajorgrids
]
\addplot [color=myorange, line width=1.2pt, forget plot]
  table[row sep=crcr]{%
0	2.57041637918865e-20\\
0.1	4.42364408458534e-20\\
0.2	7.70315521597786e-20\\
0.3	1.35729629926061e-19\\
0.4	2.40956299859524e-19\\
0.5	4.2824248049162e-19\\
0.6	7.60511613908798e-19\\
0.7	1.34784494029707e-18\\
0.8	2.38207333289671e-18\\
0.9	4.19631266822958e-18\\
1	7.36721258659262e-18\\
1.1	1.28902503408748e-17\\
1.2	2.24798071092294e-17\\
1.3	3.90822325254343e-17\\
1.4	6.77522309750088e-17\\
1.5	1.17148751777783e-16\\
1.6	2.02086897338128e-16\\
1.7	3.4788929155824e-16\\
1.8	5.97800012180608e-16\\
1.9	1.02558509624402e-15\\
2	1.75675457486602e-15\\
2.1	3.0032094339684e-15\\
2.2	5.11479534524383e-15\\
2.3	8.64497620382068e-15\\
2.4	1.44940407405978e-14\\
2.5	2.43437938422974e-14\\
2.6	4.12227317787806e-14\\
2.7	7.01671365194144e-14\\
2.8	1.19659643353919e-13\\
2.9	2.04172203652843e-13\\
3	3.48440884553709e-13\\
3.1	5.94765985930622e-13\\
3.2	1.0155389776256e-12\\
3.3	1.73477243769539e-12\\
3.4	2.96515423541666e-12\\
3.5	5.07191271020952e-12\\
3.6	8.68305957771762e-12\\
3.7	1.48799859064244e-11\\
3.8	2.55274996765182e-11\\
3.9	4.38462315994632e-11\\
4	7.54069758570561e-11\\
4.1	1.29860712434353e-10\\
4.2	2.23953304493904e-10\\
4.3	3.86786848441752e-10\\
4.4	6.69015567032715e-10\\
4.5	1.15894157623718e-09\\
4.6	2.0107262750515e-09\\
4.7	3.49390146325443e-09\\
4.8	6.08036816545271e-09\\
4.9	1.05974419920708e-08\\
5	1.84973686854044e-08\\
5.1	3.23323992527424e-08\\
5.2	5.65929193715942e-08\\
5.3	9.91871827033438e-08\\
5.4	1.74055463115499e-07\\
5.5	3.05790707471291e-07\\
5.6	5.37809057182194e-07\\
5.7	9.46801534683014e-07\\
5.8	1.66826932412905e-06\\
5.9	2.94178576902019e-06\\
6	5.19101363449063e-06\\
6.1	9.16526630767378e-06\\
6.2	1.61898605975566e-05\\
6.3	2.86086931183632e-05\\
6.4	5.05662464218138e-05\\
6.5	8.93885263204024e-05\\
6.6	0.000158019957968638\\
6.7	0.000279320307543119\\
6.8	0.000493634625175485\\
6.9	0.000872116063100823\\
7	0.00154014725480044\\
7.1	0.00271846171952012\\
7.2	0.00479523042626824\\
7.3	0.0084521866946801\\
7.4	0.0148847286668803\\
7.5	0.0261846080982656\\
7.6	0.0460016122978091\\
7.7	0.0806773782403081\\
7.8	0.141049443019349\\
7.9	0.245289231882251\\
8	0.424220530241455\\
8.1	0.728314345527349\\
8.2	1.23596865060599\\
8.3	2.05757654308903\\
8.4	3.32013770198197\\
8.5	5.10597032613159\\
8.6	7.33538263385118\\
8.7	9.67507120081411\\
8.8	11.6489865204998\\
8.9	12.9483710318986\\
9	13.5991019019695\\
9.1	13.9011432213722\\
9.2	14.1952279875499\\
9.3	14.4682125076655\\
9.4	14.717860499296\\
9.5	14.9456958267644\\
9.6	15.1540730911647\\
9.7	15.3452199241111\\
9.8	15.5210279536731\\
9.9	15.6830732497284\\
10	15.8326811473664\\
10.1	15.9709839670858\\
10.2	16.0989711668076\\
10.3	16.2175281621024\\
10.4	16.3274543899712\\
10.5	16.4294741491139\\
10.6	16.5242490334809\\
10.7	16.6123872909476\\
10.8	16.6944512640192\\
10.9	16.770963617022\\
11	16.8424127761842\\
11.1	16.9092578486425\\
11.2	16.9719306124293\\
11.3	17.0240029564529\\
11.4	17.0678614496009\\
11.5	17.1052895076502\\
11.6	17.1376223610863\\
11.7	17.1658650112462\\
11.8	17.1907787919538\\
11.9	17.21294489949\\
12	17.2328110310503\\
};
\addplot [color=mygreen, line width=1.2pt, forget plot]
  table[row sep=crcr]{%
0	5\\
0.1	4.919999899901\\
0.2	4.77999959960596\\
0.3	4.57999909912462\\
0.4	4.32548206055586\\
0.5	4.0337615993854\\
0.6	3.72475583319397\\
0.7	3.4119623514579\\
0.8	3.10428924963128\\
0.9	2.80738286956876\\
1	2.52458866259568\\
1.1	2.25764479576591\\
1.2	2.00718124098293\\
1.3	1.77307744846678\\
1.4	1.5547173577236\\
1.5	1.35117001553223\\
1.6	1.16131641051118\\
1.7	0.983937537239115\\
1.8	0.817774615870405\\
1.9	0.661569409364034\\
2	0.514090403574954\\
2.1	0.374149029372973\\
2.2	0.240608951309724\\
2.3	0.112646214502263\\
2.4	0.0103395081739279\\
2.5	0.129074176916333\\
2.6	0.244383159784235\\
2.7	0.357140997707134\\
2.8	0.468237149602402\\
2.9	0.578552871489009\\
3	0.688945721041179\\
3.1	0.800239164561063\\
3.2	0.913215492750666\\
3.3	1.02861079426602\\
3.4	1.14711114148778\\
3.5	1.26934944644714\\
3.6	1.39590267119949\\
3.7	1.52728924360699\\
3.8	1.66396664894583\\
3.9	1.80632924915764\\
4	1.95470643191982\\
4.1	2.10936121661248\\
4.2	2.27048944836483\\
4.3	2.43821969871065\\
4.4	2.61261396565465\\
4.5	2.79366923069195\\
4.6	2.98131988910481\\
4.7	3.1754410263339\\
4.8	3.37585247104498\\
4.9	3.582323518127\\
5	3.79457818519833\\
5.1	4.01230084637521\\
5.2	4.23514207814068\\
5.3	4.4627245540831\\
5.4	4.69464883694155\\
5.5	4.93049893586618\\
5.6	5.16984752162968\\
5.7	5.41226072013274\\
5.8	5.6573324111436\\
5.9	5.90467997045909\\
6	6.1539091032252\\
6.1	6.40461827897219\\
6.2	6.65640467637246\\
6.3	6.90886885876114\\
6.4	7.16161845312764\\
6.5	7.41427104217483\\
6.6	7.66645643130009\\
6.7	7.9178184158342\\
6.8	8.16801614573249\\
6.9	8.41672516309197\\
7	8.66363817088894\\
7.1	8.90846557826734\\
7.2	9.15093585765516\\
7.3	9.39079574126765\\
7.4	9.62781027891831\\
7.5	9.861762774795\\
7.6	10.0924536883855\\
7.7	10.3196988766768\\
7.8	10.5434095571473\\
7.9	10.7655415234281\\
8	10.9898046816336\\
8.1	11.218953851167\\
8.2	11.4544595264177\\
8.3	11.696953383089\\
8.4	11.9466074107686\\
8.5	12.2033934421203\\
8.6	12.4672632623159\\
8.7	12.7382837787636\\
8.8	13.0167529056832\\
8.9	13.3033161900933\\
9	13.5991015339677\\
9.1	13.9011428613692\\
9.2	14.195227635043\\
9.3	14.468212161864\\
9.4	14.7178601594127\\
9.5	14.9456954921005\\
9.6	15.1540727611306\\
9.7	15.3452195980458\\
9.8	15.5210276314918\\
9.9	15.6830729309045\\
10	15.8326808315884\\
10.1	15.9709836540756\\
10.2	16.0989708563248\\
10.3	16.2175278539331\\
10.4	16.3274540839262\\
10.5	16.429473845026\\
10.6	16.5242487312032\\
10.7	16.6123869903525\\
10.8	16.6944509645016\\
10.9	16.7709633194881\\
11	16.8424124794018\\
11.1	16.9092575539399\\
11.2	16.9719329981805\\
11.3	17.0311847866807\\
11.4	17.0879374882379\\
11.5	17.1428661545053\\
11.6	17.1964634926885\\
11.7	17.2490891202557\\
11.8	17.3010055150416\\
11.9	17.3524042545412\\
12	17.4034251674619\\
};
\draw[fill=mycolor1, draw=none, opacity=0.1] (axis cs:0,-0.4) rectangle (axis cs:2.85,18);
\draw[fill=mycolor2, draw=none, opacity=0.1] (axis cs:2.85,-0.4) rectangle (axis cs:9.1,18);
\draw[fill=mycolor3, draw=none, opacity=0.1] (axis cs:9.1,-0.4) rectangle (axis cs:12,18);

\node[right, align=left, font=\color{mygreen}] at (5,164) {\small\textbf{Actual dist.}, cf. \eqref{eq:NMPC_CA_safeDistAbs}};
\node[right, align=left, font=\color{myorange}] at (5,144) {\small\textbf{Min.~~ dist.}, cf. \eqref{eq:NMPCwCBF_safeDist_safeDist}};

\node[right, align=left, font=\color{black}] at (5,93) {\small\textbf{(1)}};
\node[right, align=left, font=\color{black}] at (53,93) {\small\textbf{(2)}};
\node[right, align=left, font=\color{black}] at (97,93) {\small\textbf{(3)}};
\end{axis}

\end{tikzpicture}%
		\vspace*{-4mm}
		\caption{Validation Scenario: The actual distance (green) is always greater than or equal to the min. safety distance (orange).} 
        \vspace*{-2mm}
		\label{fig:sim_results_overtakingDist}
\end{figure}
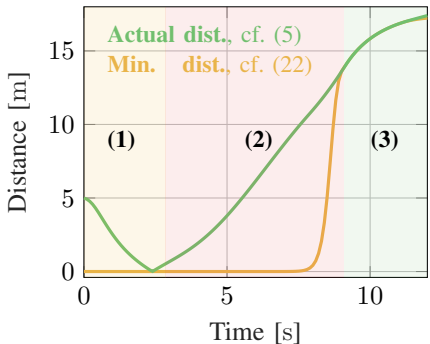

\subsection{Simulation Setup}
\label{sec:results_setup}

We investigate a realistic lane merging scenario as depicted in \prettyref{fig:problem_description_sketch} to \textsc{i}) validate our methodology 
(see \prettyref{sec:results_valScenario}) and \textsc{ii}) convey that terminal (q)DTCBF certificates may reduce the cumulative MPC cost, 
and qDTCBF certificates allow for tighter input bounds 
(see \prettyref{sec:results_costReduction}). In both cases, we use 
the following parameters: \mbox{$d_{0} = \unit[5]{m}$} (stopping distance), \mbox{$t_{h} = \unit[1]{s}$} (headway time), \mbox{$m_{lf}=10$} (leader/follower activation function),  \mbox{$\gamma_v=0.8$} (velocity DTCBF). 
Further, all agents have dimensions \mbox{$L_i=\unit[4.2]{m}$} and $W_i=\unit[2.0]{m}$.
Simulations are run on an Intel i7-7700HQ machine at \unit[3.8]{GHz} with Matlab 2022a. OCP \eqref{eq:NMPCwCBF_OCP_defOCP} is solved every $T_s=\unit[100]{ms}$ using Casadi v.3.5.5 with \texttt{ipopt} as NLP solver \cite{Andersson2019a}. To satisfy \prettyref{ass:DTCBF_background_modelUncertainties},
the simulation model and the control-oriented NMPC model 
are identical.

\begin{figure*}[ht!]
		\centering
		\setlength\fwidth{0.80\textwidth} 
		\input{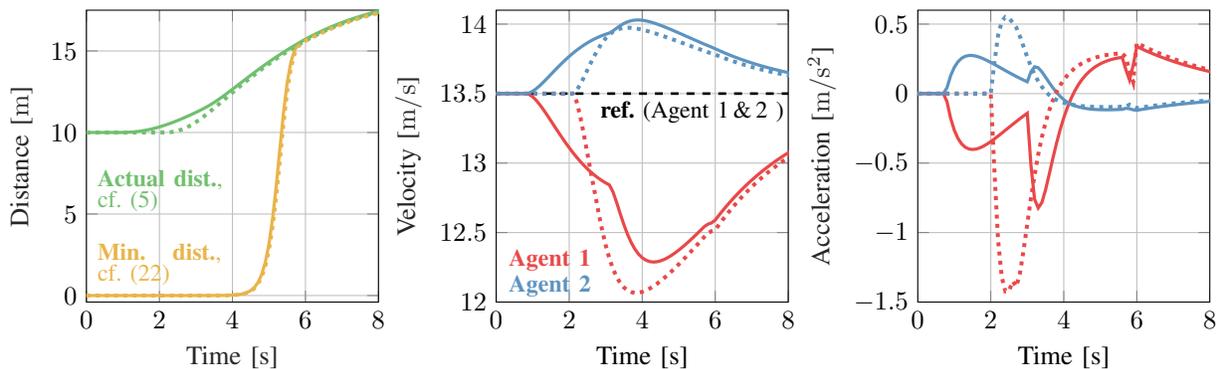}	
		\vspace*{-4.5mm}
		\caption{Simulation results for $\gamma_d=0.05$ (solid) and $\gamma_d=0.6$ (dotted) for a horizon length of $N=4$. A lower value of $\gamma_d$ results in an earlier reaction of the agents and thus in a lower input magnitude.} 
		\vspace*{-5mm}
		\label{fig:sim_results_gammaComparison}
\end{figure*}


\subsection{Validation Scenario: Overtaking and Merging in Front}
\label{sec:results_valScenario}

The first scenario demonstrates that Agent 1, which is initially behind Agent 2 (with $s_1(t_0) = \unit[-165]{m}$, \mbox{$s_2(t_0) = \unit[-160]{m}$}), is able to overtake 
and merge in front of \mbox{Agent 2} --- enabled by constraint \eqref{eq:NMPCwCBF_OCP_defOCP_safeDistNminus2}. \mbox{Agent 1} exhibits an initial velocity of $\unit[13]{m/s}$, which is slightly higher than Agent 2's initial velocity of $\unit[12.5]{m/s}$. We use constant velocity references $v_{1,\text{ref}}=\unit[13]{m/s}$ and $v_{2,\text{ref}}=\unit[12.5]{m/s}$, and choose
further parameters as: $m_{d,0} = 0.4$, \mbox{$c_{d,0} = -45$,} \mbox{$m_{d,N} = 0.06$,} \mbox{$c_{d,N} = -75$,} $\epsilon_d = 0.0025$ (safety activation function), \mbox{$\Delta \underline{v} = \unit[0.01]{m/s}$,} \mbox{$\mathbb{U}_1 = \mathbb{U}_2=[\unit[-3],\, \unit[3]]\,\unit{m/s^2}$,}  \mbox{$\overline{v} = \unit[15]{m/s}$,} 
$\gamma_d=0.15$ (safety distance qDTCBF), \mbox{$R=\mathrm{diag}(1,1)$,} $Q=Q_N=\mathrm{diag}(0,10,0,10)$ (weights), $N=15$ (horizon length). 

\prettyref{fig:sim_results_overtakingPosVelAccel} illustrates (from left to right) the position $s_i$, velocity $v_i$ and acceleration $u_i=a_i$ of Agent 1 (red) and \mbox{Agent 2} (blue). The maneuver can be subdivided into 3 phases. In phase 1, Agent 1 overtakes Agent 2 which causes Agent 1 to accelerate and Agent 2 to decelerate. Phase 2 requires the agents to accommodate an increasing minimum safety distance over the horizon, especially at the terminal stage. Finally, in phase 3 the agents proceed with velocity reference tracking, where Agent 1 is moving away from Agent 2. \prettyref{fig:sim_results_overtakingDist} depicts the actual distance (green; cf. \eqref{eq:NMPC_CA_safeDistAbs}) and the minimum safety distance (orange; cf. \eqref{eq:NMPCwCBF_safeDist_safeDist}) at current time $t_k$. Evidently, the maneuver is safe as the minimum safety distance is ensured at all times. 
Moreover, velocity and input constraints are satisfied (see \prettyref{fig:sim_results_overtakingPosVelAccel}). We can conclude that our methodology solves \prettyref{prob:problem_description_problemDef}. 

%

\subsection{Cumulative Cost and Input Bounds} 
\label{sec:results_costReduction}

The second scenario conveys that terminal (q)DTCBF certificates introduce an additional degree of freedom to reduce the cumulative MPC cost, especially on short horizons. 
Further, we demonstrate that qDTCBF certificates can be applied with lower input bounds $\lvert \underline{u} \rvert$, $\lvert \overline{u} \rvert$ compared to DTCBF certificates. We focus on safety \mbox{qDTCBF certificate} \eqref{eq:NMPCwCBF_safeDist_distCBFCons} and vary $\gamma_d$ while all other parameters remain constant, mainly because velocity DTCBF certificates have a minor influence on performance, see \prettyref{sec:NMPC_agentObjCons}. 
As initial conditions, we have chosen $s_1(t_0) = \unit[-115]{m}$, $s_2(t_0) = \unit[-105]{m}$, $v_1(t_0) = v_2(t_0) = \unit[13.5]{m}$ such that the agents have an initial distance of $\unit[10]{m}$ while driving at the same velocity. Again, we use constant velocity references  \mbox{$v_{1,\text{ref}}=v_{2,\text{ref}}=\unit[13.5]{m/s}$}, and choose further parameters as:
$m_{d,0} = 0.4$, $c_{d,0} = -45$, $m_{d,N} = 0.045$, $c_{d,N} = -85$, $\epsilon_d = 0.0025$ (safety activation function), $\mathbb{U}_1 = \mathbb{U}_2 = [\unit[-4.8],\, \unit[4.8]]\,\unit{m/s^2}$, $\overline{v} = \unit[14.5]{m/s}$, $\Delta \underline{v} = \unit[0.01]{m/s}$, $R=\mathrm{diag}(1,1)$, $Q=Q_N=\mathrm{diag}(1,0,1,0)$ (weights), $N=4$ (horizon length). 

\begin{table}[b!]
	\vspace*{-4mm}
	\begin{center}
		\caption{NMPC Cost for Different $(\gamma_d,N)$ Combinations}
  \vspace*{-1mm}
		\label{tab:results_gammaComparison}
		\begin{tabular}{ccccc}
			\hline\hline		
			 &  & Tracking Cost & Actuation Cost & Stage Cost\\
			$\gamma_d$ & $N$ & $\sum_k e_k^{\T}Q e_k$ & $\sum_k u_k^{\T}R u_k$ & $\sum_k \ell(x_k,u_k)$ \\
           \hline
             0.05 & 4 & 56.7 (-19.9\%)& \,~9.2 (-56.6\%) & 65.9 (-28.4\%)\\
             0.2 & 4 & 67.7 (~-4.4\%) & 16.7 (-21.2\%)   & 84.4 (~-8.3\%)\\
             0.4 & 4 & 69.9 (~-1.3\%) & 19.7 (~-7.1\%)   & 89.6 (~-2.6\%)\\
             0.6 & 4 & 70.8 (~~0.0\%) & 21.2 (~~0.0\%)   & 92.0 (~~0.0\%)\\
             \hline
             0.05 & 6 & 54.3 (-14.9\%) & \,~8.6 (-47.2\%)& 62.9 (-21.5\%)\\
            0.2 & 6 & 61.8 (~-3.1\%)   & 13.3 (-18.4\%)  & 75.1 (~-6.2\%)\\
             0.4 & 6 & 63.3 (~-0.8\%)  & 15.3 (~-6.1\%)  & 78.6 (~-1.9\%)\\
             0.6 & 6 & 63.8 (~~0.0\%)  & 16.3 (~~0.0\%)  & 80.1 (~~0.0\%)\\
            \hline\hline
		\end{tabular}
	\end{center}
	\vspace*{-5mm}
\end{table}

\prettyref{fig:sim_results_gammaComparison} illustrates (from left to right) the actual distance (green; cf. \eqref{eq:NMPC_CA_safeDistAbs}) and the minimum safety distance (orange; cf. \eqref{eq:NMPCwCBF_safeDist_safeDist}) at current time $t_k$, the velocity $v_i$ and the acceleration $u_i=a_i$ of Agent 1 (red) and Agent 2 (blue) for $\gamma_d=0.05$ (solid) and $\gamma_d=0.6$ (dotted). Evidently, a reduced value of $\gamma_d$ results in an earlier reaction of the control system. In fact, both agents react $\unit[1]{s}$ earlier when applying $\gamma_d=0.05$, which results in a noticeably reduced input magnitude. 

\begin{figure}[t!]
            \hspace*{5mm}
		\vspace*{0mm}
		\setlength\fwidth{0.23\textwidth}	
%
%
\definecolor{mycolor1}{rgb}{0.00000,0.44700,0.74100}%
\begin{tikzpicture}

\begin{axis}[%
width=0.958\fwidth,
height=0.75\fwidth,
at={(0\fwidth,0\fwidth)},
scale only axis,
xmin=0,
xmax=1,
xlabel style={font=\color{white!15!black}},
xlabel={$\gamma_d$ $[-]$},
ymin=1,
ymax=8.2,
ylabel style={font=\color{white!15!black}},
ylabel={$\lvert \underline{u}_i \rvert$, $\lvert \overline{u}_i \rvert$ $[\unit{m/s^2}]$},
axis background/.style={fill=white},
xmajorgrids,
ymajorgrids
]

\addplot [color=black, dashed, line width=1.2pt, forget plot]
  table[row sep=crcr]{%
0	4.8\\
0.6	4.8\\
};
\addplot [color=black, dashed, line width=1.2pt, forget plot]
  table[row sep=crcr]{%
0.6	0\\
0.6	4.8\\
};
\addplot [color=myblue, line width=1.2pt, mark=*, mark size=1.2pt, mark options={solid, myblue}, forget plot]
  table[row sep=crcr]{%
0.05	1.4\\
0.1	2.05\\
0.15	2.56\\
0.2	2.99\\
0.25	3.33\\
0.3	3.62\\
0.35	3.85\\
0.4	4.15\\
0.45	4.32\\
0.5	4.46\\
0.55	4.61\\
0.6	4.79\\
0.65	4.96\\
0.7	5\\
0.75	5.1\\
0.8	5.23\\
0.85	5.26\\
0.9	5.38\\
0.95	5.51\\
1	5.58\\
};

\addplot [color=myred, line width=1.2pt, mark=*, mark size=1.2pt, mark options={solid, myred}, forget plot]
  table[row sep=crcr]{%
0.05	7.90\\
0.1	    7.89\\
0.15	7.89\\
0.2	    7.89\\
0.25	7.89\\
0.3	    7.89\\
0.35	7.89\\
0.4	    7.89\\
0.45	7.89\\
0.5	    7.89\\
0.55	7.89\\
0.6	    7.88\\
0.65	7.88\\
0.7	    7.88\\
0.75	7.88\\
0.8	    7.88\\
0.85	7.88\\
0.9	    7.88\\
0.95	7.88\\
1	    7.88\\
};

\end{axis}

\node[right, align=left, font=\color{myred}] at (2.6,1.03) {\small\textbf{DTCBF}};
\node[right, align=left, font=\color{myblue}] at (2.42,0.65) {\small\textbf{qDTCBF}};

\end{tikzpicture}%
		\vspace*{-3.5mm}
		\caption{Least absolute input (acceleration) bounds for recursive feasibility in dependence of $\gamma_d$ -- for \prettyref{def:DTCBF_background_DTCBF} (red) and \prettyref{def:DTCBF_terminalSetTwoStep_twoStepCond} (blue). 
        Only $\gamma_d$ is varied, the other parameters (cf. \prettyref{sec:results_costReduction}) are fixed.
  With \mbox{$\mathbb{U}_i=[-4.8,\,4.8]\,\unit{m/s^2}$}, OCP \eqref{eq:NMPCwCBF_OCP_defOCP} is recursively feasible for \mbox{$\gamma_d \leq 0.6$.}} 
		\vspace*{-6mm}
		\label{fig:sim_results_aabsvsgamma}
\end{figure}
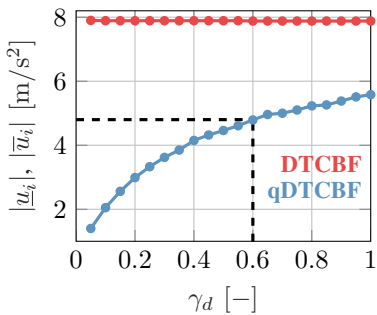

To analyze this observation more thoroughly, we have determined the cumulative tracking and actuation cost, i.e., the sum of the tracking cost $e_k^{\T}Q e_k$ (with tracking error $e_k$) and actuation cost $u_k^{\T}Ru_k$, respectively, at every time $t_k$ over the entire simulation time $0 \leq t_k \leq t_f$, where $t_f \in \mathbb{R}_{>0}$ is the duration of the simulation. The sum of the two is the overall stage cost. \prettyref{tab:results_gammaComparison} conveys the results for the horizon lengths $N=4$ and $N=6$ when varying the parameter $\gamma_d$ between $0.05$ and $0.6$. For the given input bounds, recursive feasibility can be guaranteed for $\gamma_d \leq 0.6$ (maximum value in \prettyref{tab:results_gammaComparison}), see also \prettyref{fig:sim_results_aabsvsgamma}. 

By applying a horizon length of $N=4$, the actuation cost can be reduced by $56.6\%$ for $\gamma_d=0.05$ compared to $\gamma_d=0.6$. For $N=6$, the reduction is still $47.2\%$. Thus, a lower value of $\gamma_d$ gives a lower cost but as the horizon increases the cost reduction is less significant in magnitude. Interestingly, also the tracking cost can be decreased by $19.9\%$ for \mbox{$\gamma_d=0.05$} compared to $\gamma_d=0.6$ ($N=4$) and by $14.9\%$ for $N=6$. 

\prettyref{fig:sim_results_aabsvsgamma} conveys an additional important observation, i.e., synthesizing $h_d$ as a qDTCBF by satisfying \prettyref{def:DTCBF_terminalSetTwoStep_twoStepCond}, the least absolute acceleration $\lvert \underline{u}_i \rvert$ = $\lvert \overline{u}_i \rvert$, $i\in \mathcal{A}$ 
to guarantee recursive feasibility decreases when $\gamma_d$ is decreased. When synthesizing $h_d$ as a DTCBF in accordance with \prettyref{def:DTCBF_background_DTCBF}, though, these input bounds remain essentially constant when varying $\gamma_d$. Hence, qDTCBF certificates are beneficial when only limited actuation is possible. 

%

\section{CONCLUSIONS AND FUTURE WORK}
\label{sec:conclusion}

We have proposed to utilize (q)DTCBF certificates to guarantee recursive feasibility in nonlinear MPC, and with even tighter input bounds when using qDTCBFs. By applying this methodology to lane merging, we have provided evidence that such terminal certificates offer an additional degree of freedom to reduce the cumulative MPC cost, especially on short horizons. 
Within future work, we will devise an optimization-based framework to automate the verification and synthesis of (q)DTCBFs to simplify their application.


\bibliographystyle{IEEEtran}
\bibliography{IEEEabrv}

\end{document}